# ФРИЗЫ

Алексей Панов[a], Дмитрий Панов[b], Петр Панов[a]
[a]HSE University, [b]King's College London

Фриз – это периодический рисунок на плоскости, заполняющий полосу между двумя параллельными прямыми.

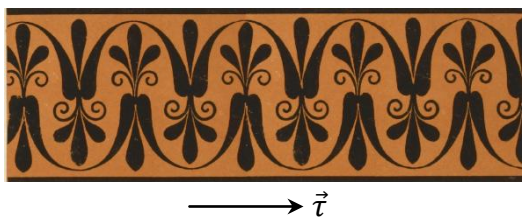

Рис. 1. Греческий фриз

На рисунке 1 представлен фрагмент фриза, на самом деле простирающего вправо и влево до бесконечности. А снизу дорисован вектор минимальной длины, при сдвиге на который фриз совмещается сам с собой. В зависимости обстоятельств, мы будем называть периодом фриза либо сам этот вектор $\vec{\tau}$, либо его длину – положительное число $\tau = |\vec{\tau}|$.

В 1856 году английский архитектор и дизайнер Оуэн Джонс (1809–1874) опубликовал свою знаменитую *Грамматику орнаментов*. Это одновременно и энциклопедия, и справочник и учебник для студентов дизайнеров, содержащий тысячи образцов фризов и других орнаментов. В отдельных главах *Грамматики* представлены фризы разных эпох и разных народов. Так что это еще и своеобразная историческая и географическая классификация фризов. Но вот у математиков существует своя собственная классификация, основанная на анализе симметрий фризов. Она позволяет разделить все фризы на семь различных типов.

Вся наша статья разбита на несколько частей. Сначала мы поговорим о движениях плоскости и симметриях геометрических объектов, расположенных на ней. Потом речь пойдет о симметриях плоской полосы, являющейся носителем фриза, и мы составим таблицу умножения ее симметрий. Далее следует основной раздел статьи, который как раз и содержит классификацию фризов. А под конец мы еще поговорим о фризах и вазах и добавим кое-что о строении химических молекул. Но начнем с еще одного примера фриза.

**Фриз-синусоида**

На самом деле, фризы в больших количествах встречаются в каждом учебнике математики – это графики периодических, в том числе, тригонометрических функций. Напомним, что функция $f(x)$, называется *периодической*, если существует такое минимальное положительное число $\tau$, что для каждого $x$ будет $f(x + \tau) = f(x)$. Число $\tau$ называется *периодом* функции $f$. При этом график периодической функции $f$ совмещается сам с собой при сдвиге на горизонтальный вектор $\vec{\tau} = (\tau, 0)$.

Возьмем, к примеру, функцию $y = \sin x$. Она периодическая с периодом $\tau = 2\pi$, ее график – синусоида. Она совмещается с собой при сдвиге на горизонтальный вектор $\vec{\tau} = (2\pi, 0)$. Чтобы можно было считать такую синусоиду фризом, нужно только указать полосу, которая эту синусоиду содержит. Поскольку наименьшее значение синуса $-1$, а наибольшее 1, естественно считать, что такая полоса ограничена прямыми $y = -1$ и $y = 1$ (рис. 2).

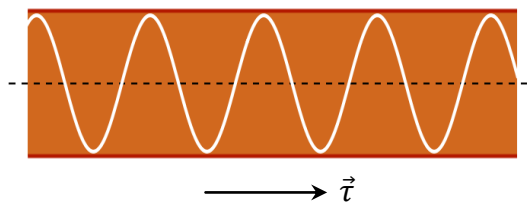

Рис. 2. Фриз-синусоида, его период равен $2\pi$

На рисунке 2 штрихом обозначена *горизонтальная ось* фриза, которая делит полосу пополам.

## ДВИЖЕНИЯ ПЛОСКОСТИ И СИММЕТРИИ

**Движения**

Движение плоскости – это преобразование плоскости, сохраняющее расстояние между точками. Об одном движении плоскости мы уже упоминали – это

- *сдвиг* плоскости на некоторый вектор $\vec{t}$, при котором каждая точка $x$ переходит в точку $x' = x + \vec{t}$

Еще одно, хорошо известное, движение плоскости – это

- *поворот* плоскости вокруг точки $A$ на угол $\varphi$

И еще одно движение – это

- *отражение относительно прямой $l$, в этом случае каждая точка $x$ переходит в точку $x'$ такую, что $x$ и $x'$ лежат по разные стороны от прямой $l$, на равных расстояниях от нее и на общем перпендикуляре к ней*

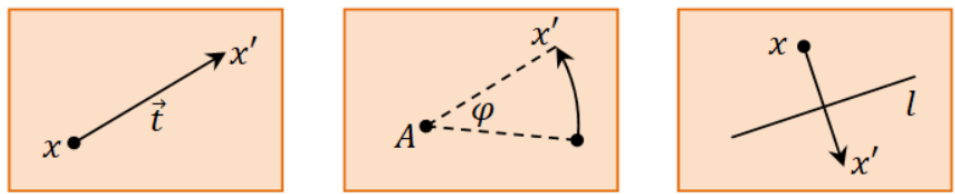

Рис. 3. Сдвиг плоскости на вектор $\vec{t}$; поворот плоскости вокруг точки $A$ на угол $\varphi$; отражение плоскости относительно прямой $l$

Наконец, добавим сюда

- *скользящее отражение* относительно прямой, которое получается в результате *отражения* относительно прямой $l$ и последующего *сдвига* на параллельный ей вектор $\vec{t}$ (рис. 4).

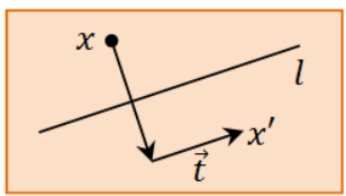

Рис. 4. Скользящее отражение плоскости,
вектор $\vec{t}$ параллелен прямой $l$

**Симметрии**

*Симметрия* объекта, расположенного на плоскости, – это движение плоскости, при котором этот объект совмещается сам с собой. Объектом, о котором идет речь в этом определении, может служить любое подмножество плоскости, например, линия, круг или конечный набор точек. Точно так же можно говорить о симметрии любой плоской картинки, в том числе, о симметрии фриза.

Упражнение 1. Перечислите все симметрии круга с центром в точке $A$.

Вот, что можно сказать по этому поводу. У круга имеется очень много симметрий, то есть существует много движений плоскости, которые переводят круг в себя. Во-первых, это все повороты плоскости вокруг центра круга, точки $A$. Во-вторых, это все отражения, относительно прямых, проходящих через центр круга. Других симметрий у круга нет.

А вот еще два упражнения, уже для самостоятельного решения.

Упражнение 2. Перечислите все симметрии а) равнобедренного треугольника, б) равностороннего треугольника.

Упражнение 3. Перечислите все симметрии а) прямоугольника, б) квадрата.

**Умножение движений и умножение симметрий**

Движения плоскости можно перемножать и их произведение тоже будет движением плоскости. Пусть, например, $R$ – это поворот плоскости на угол $\varphi$ вокруг точки $A$, а $T$ – сдвиг плоскости на вектор $\vec{t}$. Тогда их произведение $R \circ T$ – это движение, получающееся в результате последовательного применения к плоскости сначала сдвига $T$, а потом поворота $R$.

И, вообще, если имеются два движения плоскости $P$ и $Q$, то их произведение $P \circ Q$ – это движение, получающееся в результате последовательного применения к плоскости сначала движения $Q$, а потом движения $P$. Именно в таком порядке! Произведение движений $P \circ Q$ часто называют композицией этих движений.

Если два движения плоскости $P$ и $Q$ будут симметриями одного и того же геометрического объекта, то их произведение $R \circ T$ тоже будет движением и тоже будет переводить этот объект в себя и, значит, тоже будет его симметрией.

Вот два простых упражнения на эту тему.

Упражнение 4. Пусть движение $T_1$ – это сдвиг плоскости на вектор $\vec{t}_1$, а $T_2$ – сдвиг плоскости на вектор $\vec{t}_2$. Какое движение плоскости представляет их произведение $T_1 \circ T_2$?

Упражнение 5. Пусть $R_1$ и $R_2$ это две симметрии круга – повороты на углы $\varphi_1$ и $\varphi_2$ вокруг его центра. Какую симметрию круга представляет собой их произведение $R_1 \circ R_2$?

А сейчас мы обсудим более сложную задачу, а именно, поговорим о симметриях полосы и об их произведениях.

## СИММЕТРИИ ПОЛОСЫ И ТАБЛИЦА УМНОЖЕНИЯ

### Симметрии полосы

Бесконечная полоса, заключенная между двумя параллельными прямыми, является носителем фриза – фриз расположен внутри нее и целиком заполняет ее (рис. 1, 2). Любая симметрия фриза, одновременно будет и симметрией, содержащей его полосы. Так что сначала стоит отдельно разобраться именно с симметриями полосы.

Полосу мы всегда будем располагать горизонтально (рис. 5).

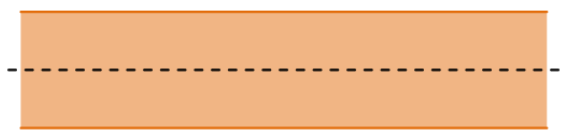

Рис. 5. Полоса со своей горизонтальной осью

Штриховая линия на рисунке – *горизонтальная ось полосы*, делящая ее пополам.

Все симметрии полосы нетрудно перечислить – это

- сдвиги в горизонтальном направлении
- повороты на 180° вокруг точек, лежащих на горизонтальной оси
- отражения относительно вертикальных прямых – вертикальных осей
- скользящие отражения относительно горизонтальной оси

Еще раз подчеркнем, что любая симметрия фриза переводит содержащую его полосу в себя и поэтому обязательно входит в этот список движений.

Введем теперь подходящие обозначения для этих симметрий полосы.

### Обозначения для симметрий полосы

Итак,

- $T_{\vec{t}}$ – *сдвиг* на вектор $\vec{t}$, параллельный горизонтальной оси (от англ. Translation)
- $R_A$ – *поворот* на угол 180° вокруг точки $A$, лежащей на горизонтальной оси (от англ. Rotation)
- $V_A$ – *вертикальное отражение,* отражение относительно вертикальной прямой, проходящей через точку $A$, лежащую на горизонтальной оси
- $S_{\vec{t}}$ – *скользящее отражение,* отражение относительно горизонтальной оси, со следующим сдвигом на вектор $\vec{t}$, параллельный горизонтальной оси (от нем. Spigel)

- $S_{\vec{0}}$ – *горизонтальное отражение*, отражение относительно горизонтальной оси, важный частный случай скользящего отражения

Введем еще одно дополнительное обозначение, которое мы дальше будем использовать

- $S'$ – это множество всех скользящих отражений относительно горизонтальной оси, за исключением горизонтального отражения $S_{\vec{0}}$

Добавим еще, что $T_{\vec{0}}$ – это тождественное отображение плоскости на себя, при котором все точки плоскости остаются на месте.

После введения этих обозначений мы можем сформулировать несколько легких упражнений на умножение симметрий полосы. В них $\vec{s}$ и $\vec{t}$ – это два горизонтальных вектора, точка $A$ лежит на горизонтальной оси полосы

<u>Упражнение 6</u>. Взгляните на определение скользящего отражения и убедитесь, что

- $S_{\vec{t}} = T_{\vec{t}} \circ S_{\vec{0}}$

<u>Упражнение 7</u>. Убедитесь также, что

- $T_{\vec{t}} \circ T_{\vec{s}} = T_{\vec{s}+\vec{t}}$
- $R_A \circ R_A = T_{\vec{0}}$ (это можно записать короче $R_A^2 = T_{\vec{0}}$)
- $S_{\vec{t}} \circ S_{\vec{s}} = T_{\vec{s}+\vec{t}}$
- $T_{\vec{t}} \circ S_{\vec{s}} = S_{\vec{t}} \circ T_{\vec{s}} = S_{\vec{s}+\vec{t}}$

Эти упражнения содержат много разных формул и, чтобы с ними было легче работать, нужно их как-то структурировать. Сейчас мы объединим их в таблицу, а именно, в таблицу умножения.

**Таблица умножения симметрий полосы**

Вот как она выглядит

Табл. 1. Таблица умножения симметрий полосы

|     | $T_{\vec{s}}$ | $R_B$ | $V_B$ | $S_{\vec{s}}$ |
|---|---|---|---|---|
| $T_{\vec{t}}$ | $T_{\vec{t}+\vec{s}}$ | $R_{B+\vec{t}/2}$ | $V_{B+\vec{t}/2}$ | $S_{\vec{t}+\vec{s}}$ |
| $R_A$ | $R_{A-\vec{s}/2}$ | $T_{2\overrightarrow{BA}}$ | $S_{2\overrightarrow{BA}}$ | $V_{A-\vec{s}/2}$ |
| $V_A$ | $V_{A-\vec{s}/2}$ | $S_{2\overrightarrow{BA}}$ | $T_{2\overrightarrow{BA}}$ | $R_{A-\vec{s}/2}$ |
| $S_{\vec{t}}$ | $S_{\vec{t}+\vec{s}}$ | $V_{B+\vec{t}/2}$ | $R_{B+\vec{t}/2}$ | $T_{\vec{t}+\vec{s}}$ |

Здесь опять принято, что точки $A$ и $B$ лежат на горизонтальной оси полосы, а векторы $\vec{s}$ и $\vec{t}$ – горизонтальны. И вот, как эта таблица работает. Предположим, что нам нужно вычислить произведение поворота $R_A$ и вертикального отражения $V_B$, то есть произведение $R_A \circ V_B$. Для этого посмотрим, что стоит в ячейке на пересечении строки с названием $R_A$ со столбцом с заголовком $V_B$. Это будет $S_{2\overrightarrow{BA}}$, и это означает, что

$$R_A \circ V_B = S_{2\overrightarrow{BA}}.$$

Упражнение 8. С помощью таблицы умножения убедитесь, что $S_{\vec{t}} \circ R_B = V_{B+\vec{t}/2}$

Конечно, остается вопрос, как были вычислены все эти 16 произведений, заполняющих эту таблицу. Но легкое Упражнение 7 гарантирует нам, что, по крайней мере, четыре клетки таблицы $T_{\vec{t}} \circ T_{\vec{s}}$, $S_{\vec{t}} \circ S_{\vec{s}}$, $T_{\vec{t}} \circ S_{\vec{s}}$, и $S_{\vec{t}} \circ T_{\vec{s}}$ заполнены правильно. А следующие два упражнения (вместе с указанием к ним) показывают, как можно проверить и всё остальное.

Упражнение 9. Докажите, что $R_A \circ V_B = S_{2\overrightarrow{BA}}$

Упражнение 10. Докажите, что $S_{\vec{t}} \circ R_B = V_{B+\vec{t}/2}$

Указание к Упражнениям 9, 10. Для доказательства этих формул воспользуйтесь следующими двумя картинками (рис. 6). На них $V_B(x)$ – это образ точки $x$ при отражении от вертикальной прямой, проходящей через точку $B$, а $R_A(V_B(x))$ – это образ точки $x$ при движении $R_A \circ V_B$. И точно так же $R_B(x)$ – это образ точки $x$ при повороте $R_B$, а $S_{\vec{t}}(R_B(x))$ – это образ точки $x$ при движении $S_{\vec{t}} \circ R_B$.

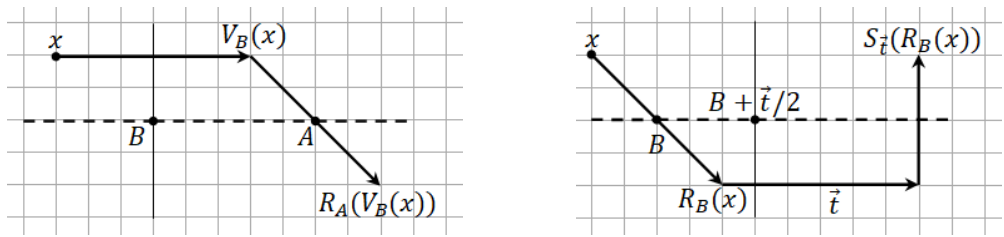

Рис. 6. К упражнениям 9, 10

И теперь заключительное упражнение на эту тему.

Упражнение 11. Докажите, что все ячейки в Таблице 1 заполнены правильно.

**Сжатая таблица умножения симметрий полосы**

Наша предыдущая таблица умножения содержит важную информацию, но выглядит достаточно громоздко. В принципе ее можно слегка сжать, убрав нижние индексы (табл. 2). И это тоже ценная конструкция. Покажем, как читать такую таблицу умножения. Например, что значит та же самая формула, содержащаяся в этой таблице,

$$R \circ V = S.$$

Она означает, что какой бы поворот $R$ и какое бы вертикальное отражение $V$, являющиеся симметриями полосы, мы ни взяли, в результате умножения $R$ на $V$ получим новую симметрию полосы, и это обязательно будет некоторое скользящее отражение $S$.

Табл. 2. Сжатая таблица умножения

|   | $T$ | $R$ | $V$ | $S$ |
|---|---|---|---|---|
| $T$ | $T$ | $R$ | $V$ | $S$ |
| $R$ | $R$ | $T$ | $S$ | $V$ |
| $V$ | $V$ | $S$ | $T$ | $R$ |
| $S$ | $S$ | $V$ | $R$ | $T$ |

Теперь, после всей проделанной работы, мы готовы перейти к классификации фризов.

## СЕМЬ ТИПОВ ФРИЗОВ

Мы будем характеризовать фриз набором его симметрий. И чтобы разобраться, что это означает, начнем с фриза-синусоиды (рис. 2).

**Симметрии фриза-синусоиды**

Во-первых, у него, как и любого другого фриза, по определению, есть горизонтальный сдвиг $T_{\vec{\tau}}$ на период, в данном случае на $\vec{\tau} = (2\pi, 0)$, который переводит его в себя, то есть является его симметрией. На самом деле у фриза, имеется бесконечно много сдвигов-симметрий и все они имеют вид $T_{n\vec{\tau}}$, где $n$ любое целое.

Во-вторых, обратите внимание на какую-либо точку пересечения синусоиды с горизонтальной осью фриза. Если повернуть весь фриз вокруг этой точки на $180°$, то он совместится с собой. Так что этот поворот – это еще одна симметрия. И таких поворотов симметрий тоже бесконечно много. На рисунке 7 они отмечены синими точками.

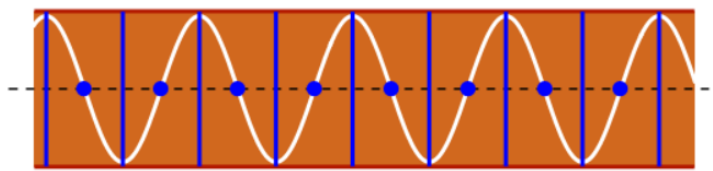

Рис. 7. Отмечены центры поворотов-симметрий,
и оси вертикальных отражений-симметрий

В-третьих, проведите вертикальные прямые через точки максимума и минимума синусоиды (рис. 7). Отражения фриза относительно этих прямых тоже являются его симметриями.

Наконец, отразите фриз-синусоиду относительно горизонтальной оси, а потом еще дополнительно сдвиньте на полпериода вправо. При таком движении фриз совместится с собой. Это скользящее отражение $S_{\vec{\tau}/2}$ тоже является симметрией этого фриза, а вместе с ним являются симметриями и все скользящие отражения вида $S_{(n+1/2)\vec{\tau}}$, где $n$ любое целое. Заметим, что среди этих скользящих отражений отсутствует горизонтальное отражение $S_{\vec{0}}$.

Нетрудно убедиться, что мы перечислили все возможные типы симметрии фриза-синусоиды: сдвиги, повороты, скользящие отражения ($S_{\vec{0}}$ среди них отсутствует), а также вертикальные отражения. Мы будем говорить, что фриз с таким набором симметрий имеет тип

$$\langle T, R, V, S' \rangle.$$

<u>Упражнение 12</u>. Убедитесь, что греческий фриз на рисунке 1 имеет тот же самый тип симметрии.

<u>Упражнение 13</u>. Убедитесь, что график-фриз функции $y = \sin 2x$ и график-фриз функции $y = \sin \frac{x}{2}$ имеют тот же самый тип симметрии, что и график-фриз функции $y = 2 \sin x$.

А сейчас мы поэтапно выясним, какие именно симметрии могут быть у фриза с заданным периодом.

**Какие сдвиги-симметрии есть у фриза**

По определению, данному в самом начале статьи, в состав симметрий любого фриза обязательно входит сдвиг $T_{\vec{\tau}}$, где вектор $\vec{\tau}$ – это период фриза, то есть горизонтальный вектор минимальной длины, при сдвиге на который фриз совмещается сам с собой. Понятно, что и любой сдвиг вида $T_{n\vec{\tau}}$ на вектор $n\vec{\tau}$, являющийся целым кратным периода, тоже будет симметрией этого фриза. Теперь объясним, почему других сдвигов-симметрий у фриза нет.

Предположим, что некоторый сдвиг $T_{\vec{t}}$ входит в состав симметрий фриза. Тогда в состав симметрий входит и любой сдвиг вида $T_{n\vec{\tau}} \circ T_{\vec{t}}$. Таблица 1 позволяет вычислить это произведение

$$T_{n\vec{\tau}} \circ T_{\vec{t}} = T_{\vec{t}+n\vec{\tau}}.$$

Если горизонтальный вектор $\vec{t}$ не кратен периоду $\vec{\tau}$, то среди векторов $\vec{t} + n\vec{\tau}$, где $n$ любое целое, обязательно найдется такой ненулевой вектор, что $|\vec{t} + n\vec{\tau}| < |\vec{\tau}|$. А это уже противоречит свойству минимальности длины периода $\vec{\tau}$. Поэтому вектор $\vec{t}$ обязательно кратен периоду $\vec{\tau}$. То есть обязательно существует такое целое $n$, что $\vec{t} = n\vec{\tau}$.

Итак, все сдвиги-симметрии фриза с периодом $\vec{\tau}$ имеют вид $T_{n\vec{\tau}}$, $n$ – целое. Других сдвигов-симметрий у фриза нет.

### Какие повороты-симметрии могут быть у фриза

Пусть у фриза с периодом $\vec{\tau}$ имеется хотя бы один поворот-симметрия. Обозначим центр этого поворота на 180°, лежащий на горизонтальной оси фриза, буквой $O$. Тогда сам поворот запишется как $R_O$. Вместе с $R_O$ симметрией фриза будет и любое произведение вида $T_{n\vec{\tau}} \circ R_O$. Из таблицы 1 извлекаем, что

$$T_{n\vec{\tau}} \circ R_O = R_{O+n\vec{\tau}/2}.$$

Таким образом, вместе с поворотом $R_O$ мы получаем целую серию поворотов симметрий $R_{O+n\vec{\tau}/2}$. Посмотрим, почему других поворотов-симметрий у фриза нет.

Пусть у фриза имеется поворот-симметрия вида $R_{O+\vec{t}}$. Тогда у него есть симметрия $R_{O+\vec{t}} \circ R_O$. Таблица 1 дает нам, что это сдвиг на вектор $2\vec{t}$

$$R_{O+\vec{t}} \circ R_O = T_{2\vec{t}}.$$

Из предыдущего пункта мы знаем, что этот вектор $2\vec{t}$ обязательно кратен периоду $\vec{\tau}$, то есть $2\vec{t} = n\vec{\tau}$ и $\vec{t} = n\vec{\tau}/2$.

Итак, фриз с периодом $\vec{\tau}$ может иметь повороты-симметрии вида $R_{O+n\vec{\tau}/2}$, где $n$ – целое, а точка $O$ – это какая-то на горизонтальной оси фриза. Других поворотов-симметрий у фриза быть не может.

### Какие вертикальные отражения-симметрии могут быть у фриза

Рассмотрим фриз, в состав симметрий которого, кроме сдвигов $T_{n\vec{\tau}}$, входят вертикальные отражения. И пусть $V_A$ – это одно из них. Обращаясь к таблице 1, находим, что

$$T_{n\vec{\tau}} \circ V_A = V_{A+n\vec{\tau}/2},$$

и что для любого горизонтального вектора $\vec{t}$

$$V_{A+\vec{t}} \circ V_A = T_{2\vec{t}}.$$

Но это фактически те же самые два соотношения, которые мы использовали в предыдущем разделе. Только буква $R$ в них заменена на букву $V$. Так что мы можем один к одному переписать предыдущий раздел, сделав эту замену, и еще заменив всюду слово *поворот* на словосочетание *вертикальное отражение*. В конце концов, мы получим следующее верное утверждение.

Фриз с периодом $\vec{\tau}$ может иметь вертикальные отражения-симметрии вида $V_{A+n\vec{\tau}/2}$, где $n$ – целое, а точка $A$ – это какая-либо точка на горизонтальной оси фриза. Других вертикальных отражений-симметрий у фриза быть не может.

### Какие скользящие отражения-симметрии могут быть у фриза

На этот раз обсудим фризы, в состав симметрий которых входят скользящие отражения. Пусть $S_{\vec{t}}$ – одно из них. Вместе с ним симметриями фриза будут и все скользящие отражения вида

$$S_{\vec{t}+n\vec{\tau}} = T_{n\vec{\tau}} \circ S_{\vec{t}}$$

Среди всех горизонтальных векторов $\vec{t} + n\vec{\tau}$ выберем вектор $\vec{t}'$ с минимальной длиной. Тогда обязательно

$$0 \leq |\vec{t}'| < |\vec{\tau}|.$$

С другой стороны

$$S_{\vec{t}'} \circ S_{\vec{t}'} = T_{2\vec{t}'}$$

Поэтому $T_{2\vec{t}'}$ – это тоже симметрия фриза, откуда следует, что вектор $2\vec{t}'$ кратен вектору периода $\vec{\tau}$. Вместе с неравенством $0 \leq |\vec{t}'| < |\vec{\tau}|$ это говорит о том, что либо $\vec{t}' = \vec{0}$, либо $\vec{t}'$ с точностью до знака совпадает с $\vec{\tau}/2$.

В первом случае, когда $\vec{t}' = \vec{0}$, отражение $S_{\vec{0}}$ будет симметрией фриза, и все скользящие отражения-симметрии фриза будут иметь вид $S_{n\vec{\tau}}$.

Во втором случае, когда $|\vec{t}'| = 2|\vec{\tau}|$, скользящее отражение $S_{\vec{\tau}/2}$ будет симметрией фриза, и все скользящие отражения-симметрии будут иметь вид $S_{(n+1/2)\vec{\tau}}$, и не будут включать $S_{\vec{0}}$.

После этого, слегка затянувшегося введения, мы можем приступить к перечислению типов фризов.

**Первые пять типов фризов** $\langle T \rangle, \langle T, R \rangle, \langle T, V \rangle, \langle T, S_{\vec{0}} \rangle, \langle T, S \rangle$

Всюду дальше вектор $\vec{\tau}$ – это период фриза.

1) Пусть в состав симметрий фриза входят только сдвиги. Мы уже знаем, что все такие сдвиги обязаны иметь вид $T_{n\vec{\tau}}$. И вот два примера подобных фризов

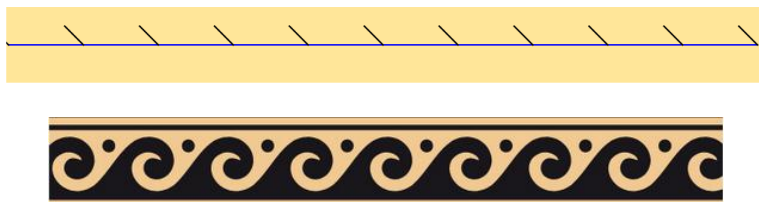

Рис. 8. Фризы типа $\langle T \rangle$

Мы говорим, что все такие фризы имеют тип $\langle T \rangle$.

2) Пусть в состав симметрий фриза входят только сдвиги и повороты, тогда, как мы знаем, все эти повороты имеют вид $R_{O+n\vec{\tau}/2}$, где $O$ – это некоторая точка на горизонтальной оси фриза. Приведем два примера таких фризов

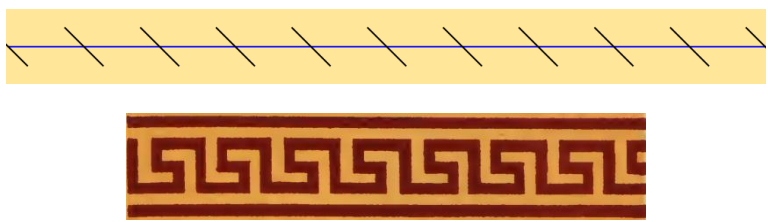

Рис. 9. Фризы типа $\langle T, R \rangle$

Все симметрии этих фризов – это сдвиги вида $T_{n\vec{\tau}}$ и повороты $R_{O+n\vec{\tau}/2}$. Мы говорим, что такие фризы имеют тип $\langle T, R \rangle$.

3) Пусть в состав симметрий фриза входят только сдвиги и вертикальные отражения. Мы знаем, что в этом случае все эти вертикальные отражения имеют вид $V_{A+n\vec{\tau}/2}$, где $A$ – это некоторая точка на горизонтальной оси фриза. Вот два таких фриза

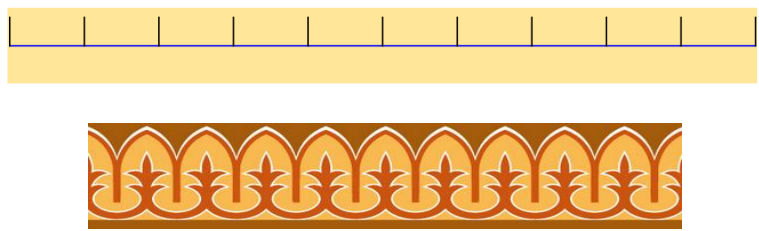

Рис. 10. Фризы типа $\langle T, V \rangle$

Все их симметрии – это сдвиги $T_{n\vec{\tau}}$ и вертикальные отражения $V_{A+n\vec{\tau}/2}$. Мы говорим, что эти фризы имеют тип $\langle T, V \rangle$.

4) Рассмотрим теперь случай, когда в состав симметрий фриза входят только сдвиги и скользящие отражения, среди которых присутствует и отражение $S_{\vec{0}}$. Раньше мы уже выяснили, что в этом случае все скользящие отражения-симметрии фриза – это $S_{n\vec{\tau}}$. И вот два примера фризов такого рода

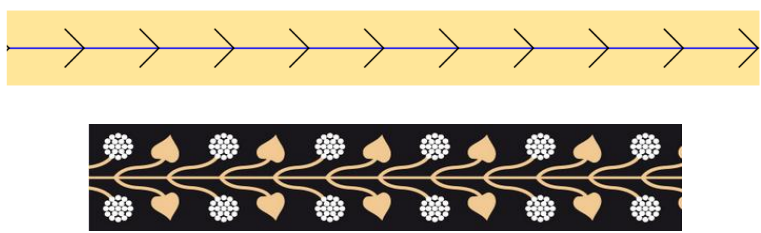

Рис. 11. Два фриза типа $\langle T, S_{\vec{0}} \rangle$

Перечислим все симметрии этих фризов. Это сдвиги $T_{n\vec{\tau}}$ и скользящие отражения $S_{n\vec{\tau}}$. Тип симметрии таких фризов мы обозначим $\langle T, S_{\vec{0}} \rangle$.

5) Наконец, посмотрим, что происходит, когда опять в состав симметрий фриза входят только сдвиги и скользящие отражения, но при этом отражение $S_{\vec{0}}$ отсутствует. На этот раз, как мы выяснили раньше, все возможные скользящие отражения-симметрии фриза – это $S_{(n+1/2)\vec{\tau}}$. И такие фризы существуют – вот два примера

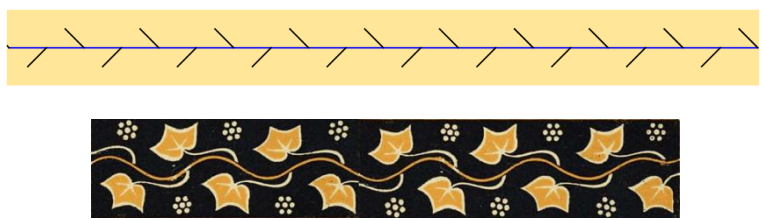

Рис. 12. Два фриза типа $\langle T, S' \rangle$

Все симметрии этих фризов – это сдвиги $T_{n\vec{\tau}}$ и скользящие отражения $S_{(n+1/2)\vec{\tau}}$. Тип таких фризов мы обозначим $\langle T, S'\rangle$.

**Еще два типа фризов – $\langle T, R, V, S_{\vec{0}}\rangle$ и $\langle T, R, V, S'\rangle$**

Из сжатой таблицы умножения симметрий (табл. 2) можно извлечь следующие три формулы

$$R \circ S = V, \quad S \circ V = R, \quad V \circ R = S.$$

Наличие этих формул позволяет сделать следующий вывод. Если у фриза в качестве симметрий присутствует какая-то пара движений из тройки вида $R, S, V$, то у него в качестве симметрии присутствуют все три типа движение из этой тройки. Так что нам осталось рассмотреть только следующие две ситуации.

6) Рассмотрим случай, когда в состав симметрий фриза с периодом $\vec{\tau}$, кроме сдвигов, входит некоторый поворот $R_O$ и отражение относительно горизонтальной оси $S_{\vec{0}}$ и, следовательно, также некоторое вертикальное отражение. Мы уже знаем, что в этом случае все сдвиги имеют вид $T_{n\vec{\tau}}$, все повороты – $R_{O+n\vec{\tau}/2}$, скользящие отражения – $S_{n\vec{\tau}}$. Что касается вертикальных отражений, то из таблицы 1 извлекаем следующее

$$R_{O+k\vec{\tau}/2} \circ S_{l\vec{\tau}} = V_{O+(k-l)\vec{\tau}/2},$$

то есть оси вертикального отражения, проходят через те же самые центры вращений.

Итак, в рассматриваемом случае все симметрии фриза – это сдвиги $T_{n\vec{\tau}}$, повороты $R_{O+n\vec{\tau}/2}$, скользящие отражения $S_{n\vec{\tau}}$ и вертикальные отражения вида $V_{O+n\vec{\tau}/2}$. При этом вертикальные оси отражения проходят через центры вращения. И вот два примера таких фризов.

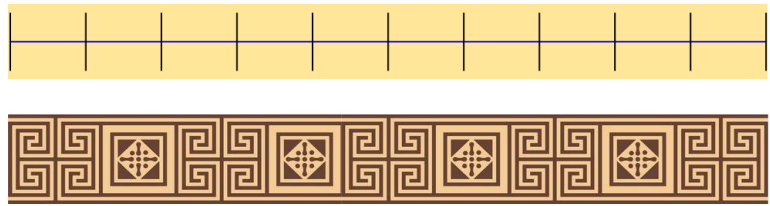

Рис. 13. Два фриза типа $\langle T, R, V, S_{\vec{0}}\rangle$

7) Наконец посмотрим, что происходит, когда в состав симметрий фриза с периодом $\vec{\tau}$, кроме сдвигов, входит некоторый поворот $R_O$, некоторое скользящее отражение, при этом отражение $S_{\vec{0}}$ отсутствует, и опять автоматически присутствует некоторое вертикальное отражение. Мы уже знаем, что в этом случае все сдвиги имеют вид $T_{n\vec{\tau}}$, все повороты – $R_{O+n\vec{\tau}/2}$, а скользящие отражения – $S_{(n+1/2)\vec{\tau}/2}$. Что касается вертикальных отражений, то с помощью таблицы 1 на этот раз получаем следующий результат

$$R_{O+k\vec{\tau}/2} \circ S_{(l+1/2)\vec{\tau}/2} = V_{(O+\vec{\tau}/4)+(k-l)\vec{\tau}/2}.$$

Итак, в рассматриваемом случае все симметрии фриза – это сдвиги $T_{n\vec{\tau}}$, повороты $R_{O+n\vec{\tau}/2}$, скользящие отражения $S_{(l+1/2)\vec{\tau}/2}$ и вертикальные отражения вида $V_{(O+\vec{\tau}/4)+n\vec{\tau}/2}$. При этом центры поворотов и вертикальные оси отражения чередуются и отстоят друг от друга на четверть периода, то есть на расстояние $\tau/4$. Вот два примера таких фризов

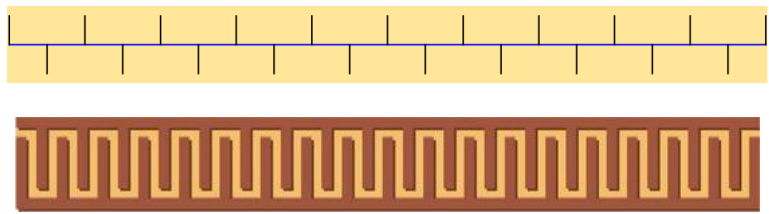

Рис. 13. Два фриза типа $\langle T, R, V, S' \rangle$

Нижний фриз, конечно, сильно напоминает фриз-синусоиду (рис. 2). К тому же самому типу симметрии относится и греческий фриз, расположенный на рисунке 1.

Таким образом, классификация фризов завершена – мы описали все семь типов их симметрий. И теперь предлагаем большое упражнение на эту тему.

Упражнение 14. Найдите в сети книгу *Owen Jones, The Grammar of Ornament* и постарайтесь определить, к какому именно типу симметрии

$$\langle T \rangle, \langle T, R \rangle, \langle T, V \rangle, \langle T, S_{\vec{0}} \rangle, \langle T, S' \rangle, \langle T, R, V, S_{\vec{0}} \rangle \text{ или } \langle T, R, V, S' \rangle$$

относится тот или иной фриз из этой книги.

Подходящим инструментом для решения этого упражнения может послужить графический редактор Paint.

**Paint – инструмент для изучения симметрий фризов**

Одно из важных для нас средств работы с изображениями в Paint называется *Повернуть*. На рисунке 14 справа оно представлено в развернутом виде. Нас, конечно же, интересуют последние три пункта. Если выделить на экране Paint расположенный там фриз, то с их помощью мы сможем повернуть этот фриз на $180°$, отразить его относительно горизонтали и отразить относительно вертикали. Что касается сдвигов выделенного фриза, то они элементарно осуществляются клавишами со стрелками.

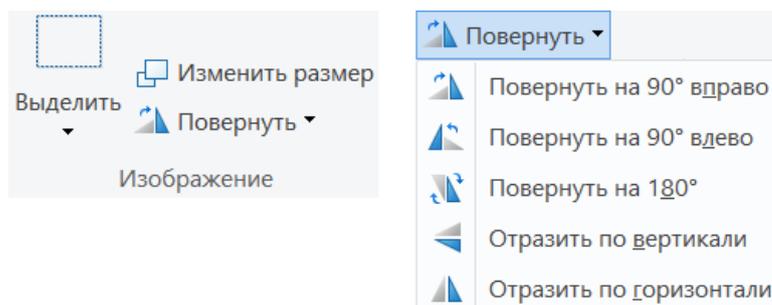

Рис. 14. Работа с изображениями в Paint

Другие интересные для нас средства Paint – это *Изменить размер* и *Наклон* (рис. 15).

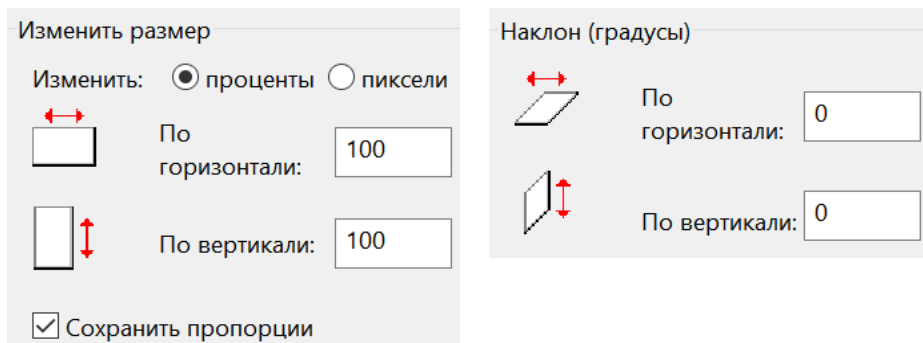

Рис. 14. Панели *Изменить размер* и *Наклон*

При наличии галочки в окошке *Сохранить пропорции*, размеры изображения по горизонтали и вертикали меняются согласованно. При снятии галочки они меняются независимо.

<u>Упражнение 14</u>. Какие из семи типов симметрий фризов меняются, какие сохраняются

- при согласованном изменении размеров по горизонтали и вертикали
- при изменении размеров только по горизонтали
- при изменении размеров только по вертикали
- при наклоне по горизонтали

**РАЗНОЕ**

**Вазы и фризы**

Греческие фризы в *Грамматике* Джонса – это не прямые копии какого-либо оригинала. В большинстве своем они срисованы с античных ваз из Британского музея и Лувра (рис. 15).

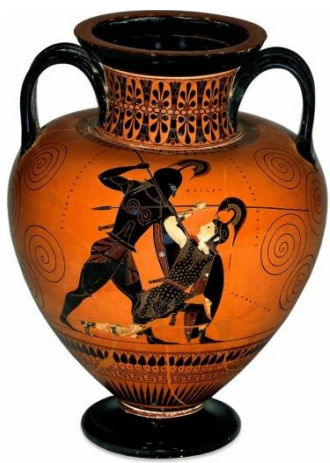

Рис. 15. Ахилл поражает Пентесилею
© The Trustees of the British Museum, CC BY-NC-SA 4.0

При этом цилиндрический узор, опоясывающий горловину вазы или расположенный вблизи ее днища, в процессе прорисовки фриза выпрямляется и периодически размножается влево и вправо, целиком заполняя горизонтальную полосу на плоскости.

И наоборот, рисунок фриза легко переносится на боковую поверхность цилиндра. Нужно только из полосы фриза вырезать кусок длины в $n$ периодов и склеить его концы.

Получится цилиндрическое кольцо с соответствующим периодическим рисунком (табл. 3, столбец 3), которое как раз и можно рассматривать как боковую поверхность цилиндра.

При этом существует естественное соответствие между симметриями фриза и симметриями порожденного таким образом цилиндра

- сдвиг фриза на период $\vec{\tau}$ → поворот цилиндра вокруг вертикальной оси на угол $2\pi/n$
- поворот фриза на 180° → поворот цилиндра на 180° вокруг соответствующей горизонтальной оси, проходящей через центр цилиндра
- скользящее отражение фриза → отражение цилиндра относительно горизонтальной плоскости симметрии и последующий поворот вокруг вертикальной оси
- вертикальное отражение фриза → отражение от соответствующей вертикальной плоскости, проходящей через ось цилиндра

Табл. 3. Фризы → цилиндры → молекулы

| Тип фриза | Фриз | Цилиндрический аналог, $n=6$ | Молекула |
|---|---|---|---|
| $\langle T \rangle$ | 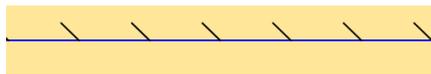 | 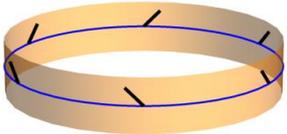 | 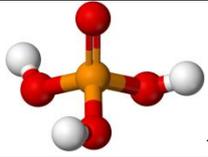 $n=3$ |
| $\langle T, R \rangle$ | 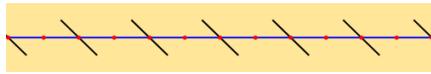 | 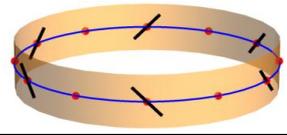 | 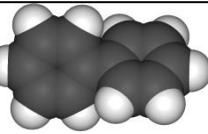 $n=2$ |
| $\langle T, V \rangle$ | 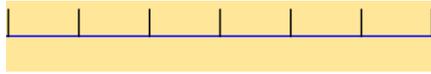 | 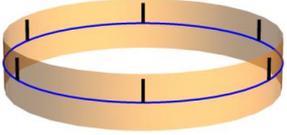 | 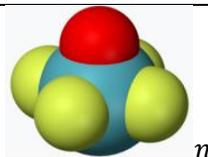 $n=4$ |
| $\langle T, S_{\vec{0}} \rangle$ | 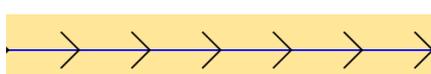 | 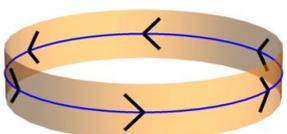 | 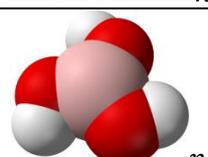 $n=3$ |
| $\langle T, S' \rangle$ | 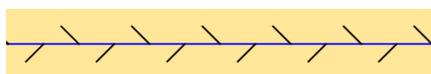 | 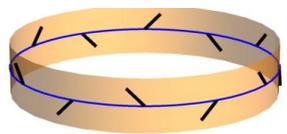 | 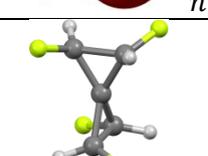 $n=2$ |
| $\langle T, R, V, S_{\vec{0}} \rangle$ | 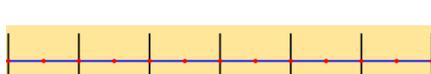 | 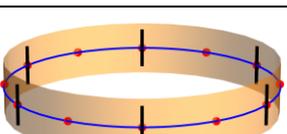 | 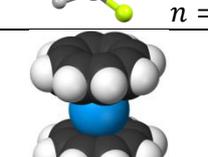 $n=8$ |
| $\langle T, R, V, S' \rangle$ | 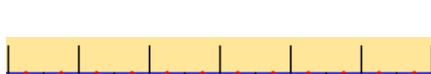 | 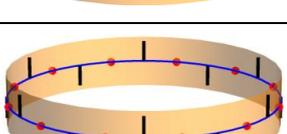 | 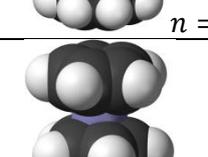 $n=5$ |

**Фризы и молекулы**

Наряду с греческими вазами существуют и другие материальные объекты, наделенные указанными цилиндрическими симметриями, например, химические молекулы (табл. 3, столбец 4). Если исключить молекулы, все атомы, которых лежат на одной прямой, то все остальные молекулы, обладающие цилиндрической симметрией, разбиваются, как и фризы, на семь типов (табл. 3). Каждая такая молекула дополнительно снабжается целым параметром $n$, определяющим наименьший положительный угол $2\pi/n$, при повороте на который молекула совмещается сама с собой.

Добавим, что наряду с молекулами, обладающими цилиндрической симметрией, существует еще один класс молекул, симметрии которых связаны с симметриями правильных многогранников. Дальнейшие интересные химические подробности вы можете почерпнуть из статьи *Молекулярная симметрия* в Википедии. Именно оттуда взяты все изображения молекул, содержащиеся в таблице 3.

Еще один ресурс, который мы рекомендуем посетить, – это веб-сайт Университета Оттербейн – *Symmetry@Otterbein*, предназначенный для студентов химиков. Это тоже хорошая возможность для знакомства с симметрией молекул и для упражнений по распознаванию их типов симметрий.

**Периодические мозаики и двумерная кристаллография**

Половина всех орнаментов, представленных в *Грамматике* Джонса, – это фризы, и есть еще одна большая составляющая – это периодические мозаики, заполняющие плоскость. В состав симметрий такой мозаики обязательно входят сдвиги *разных направлений* (рис. 14).

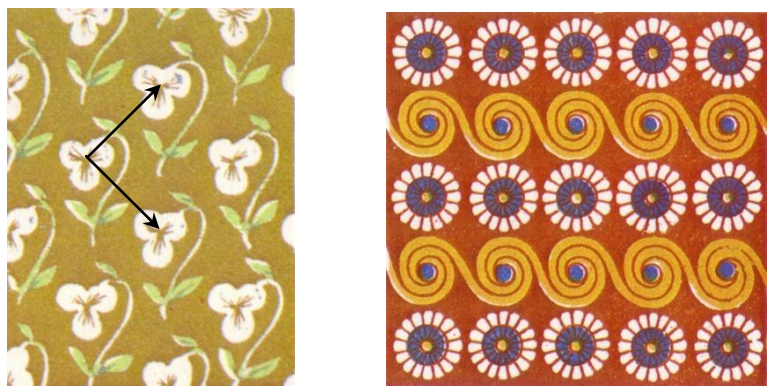

Рис. 14. Мозаики из Грамматики Джонса, для первой указаны два ее непараллельных сдвига-симметрии

И тут тоже, как для фризов, возникает естественная задача перечисления всех возможных типов симметрий таких мозаик.

В 1891 году в 28-м томе *Записок Императорского Санкт-Петербургского Минералогического Общества* была опубликована статья Евграфа Степановича Фёдорова *Симметрия правильных фигур*, в которой он практически завершил свою знаменитую классификацию трехмерных кристаллических структур. Одновременно с Фёдоровым эту зада-чу решил немецкий математик Артур Шёнфлис, и количество этих структур, подсчитанное обоими учеными, оказалось равным 230.

Для нас важно, что в том же самом 28-м томе *Записок* Фёдоров опубликовал еще одну свою статью *Симметрия на плоскости*, где перечислил и все двумерные кристаллические структуры, общее число которых оказалось равным 17.

На самом деле, задача перечисления всех возможных двумерных кристаллических структур и задача перечисления периодических мозаик на плоскости – это одна и та же задача. Для примера, посмотрите на две мозаики, размещенные на рисунке. На самом деле, – это карты функции локализации электрона в определенных кристаллографических плоскостях для систем Ga–Te при высоких давлениях[1].

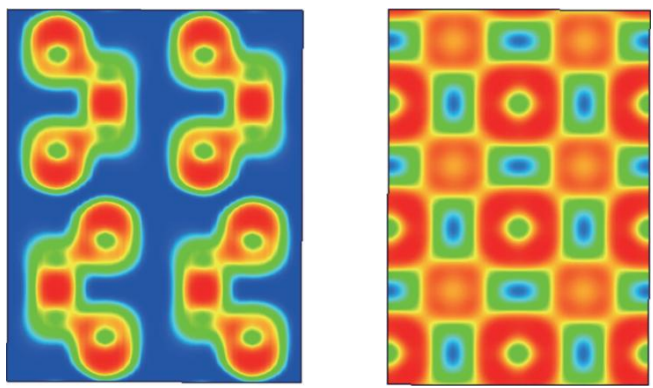

Рис. 15. Две мозаики – цветные карты функции локализации электрона,
(воспроизводится с разрешения авторов)

Так что Фёдоров одновременно с перечислением двумерных кристаллических структур, перечислил и все типы симметрий периодических мозаик на плоскости.

Кстати, о мозаиках достаточно подробно рассказано в книге Гарольда Коксетера *Введение в геометрию*, в её четвертой главе *Двумерная кристаллография*.

Добавим еще, что основным алгебраическим инструментом исследования симметрий геометрических объектов являются группы. Но мы, специально, язык теории групп здесь не использовали.

На этом наш рассказ закончен. В качестве дополнительных источников по этой теме мы рекомендуем

- Owen Jones, *The Grammar of Ornament*, 1856
- Коксетер, *Введение в геометрию*, 1966
  Глава 3  *Движения в евклидовой плоскости*,
  Глава 4  *Двумерная кристаллография*

---

[1] Youchun Wang, Fubo Tian *et al, Structural and electrical properties of Ga–Te systems under high pressure*,
2019 *Chinese Phys. B* **28** 056104

- *Фейнмановские лекции по физике*, т 7, 1977,
  Глава 30, *Внутренняя геометрия кристаллов*,
  §1 *Внутренняя геометрия кристаллов*,
  §5 *Симметрии в двух измерениях*
- Г. Вейль, *Симметрия*, 1968,
  Глава 3  *Орнаментальная симметрия*
- А.А. Ошемков, Ф.Ю. Попеленский, А.А. Тужилин, А.Т. Фоменко, А.И. Шафаревич,
  *Курс наглядной геометрии и топологии*, 2014,
  Глава 14. *Симметрии плоских кристаллов*
- А.И. Шафаревич, *Курс наглядной геометрии и топологии,* видеолекции, 2019,
  Лекция 5. *Группы движений плоскости*
  Лекция 6. *Классификация двумерных кристаллов*
  https://teach-in.ru/course/visual-geometry-and-topology-shafarevich